

\input amstex
\documentstyle{amsppt}

\loadbold
\loadeurm
\loadeurb
\loadeusm

\magnification=\magstep1
\hsize=6.7truein
\vsize=9.6truein
\hcorrection{-0.1truein} 
\vcorrection{-0.2truein} 

\font\sc=cmcsc10
  
\font \smallrm=cmr10 at 11truept
 at 7truept
\font \smallbf=cmbx10 at 11truept 
 at 11truept
\font \smallsl=cmsl10 at 11truept

\baselineskip=14pt

\def \longtwoheadrightarrow
 {\relbar\joinrel\relbar\joinrel\twoheadrightarrow}
\def \longhookrightarrow {\lhook\joinrel\relbar\joinrel\rightarrow}
\def \llonghookrightarrow
 {\lhook\joinrel\relbar\joinrel\relbar\joinrel\relbar\joinrel\rightarrow}

\def \N {\Bbb{N}}
\def \Z {\Bbb{Z}}
\def \C {\Bbb{C}}
\def \Zeps {\Z_\varepsilon}

\def \Oqg {\Cal{O}_q(G)}

\def \Oeg {\Cal{O}_\varepsilon(G)}
\def \Oezg {\Cal{O}_\varepsilon^{\,\Zeps\!}(G)}
\def \Ounog {\Cal{O}_1(G)}

\def \Og {\Cal{O}(G)}

\def \Ozeg {\Cal{O}^{\,\Zeps\!}(G)}
\def \Oqgl {\Cal{O}_q(GL_n)}
\def \Oqrgl {\Cal{O}_q^{\,R}({GL}_n)}
\def \Oqzgl {\Cal{O}_q^{\,\Bbb{Z}_q}({GL}_n)}
\def \Oegl {\Cal{O}_\varepsilon({GL}_n)}
\def \Oezgl {\Cal{O}_\varepsilon^{\,\Zeps\!}({GL}_n)}

\def \Ounozgl {\Cal{O}_1^{\,\Bbb{Z}}({GL}_n)}
\def \Ogl {\Cal{O}({GL}_n)}
\def \Ozgl {\Cal{O}^{\,\Bbb{Z}}({GL}_n)}
\def \Ozegl {\Cal{O}^{\,\Zeps\!}({GL}_n)}
\def \Oqsl {\Cal{O}_q({SL}_n)}
\def \Oqrsl {\Cal{O}_q^{\,R}({SL}_n)}
\def \Oqzsl {\Cal{O}_q^{\,\Bbb{Z}_q}({SL}_n)}

\def \Oezsl {\Cal{O}_\varepsilon^{\,\Zeps\!}({SL}_n)}

\def \Ounozsl {\Cal{O}_1^{\,\Bbb{Z}}({SL}_n)}

\def \Ozsl {\Cal{O}^{\,\Bbb{Z}}({SL}_n)}
\def \Ozesl {\Cal{O}^{\,\Zeps\!}({SL}_n)}
\def \Oqm {\Cal{O}_q(M_n)}
\def \Oqrm {\Cal{O}_q^{\,R}(M_n)}
\def \Oqzm {\Cal{O}_q^{\,\Bbb{Z}_q}(M_n)}
\def \Oem {\Cal{O}_\varepsilon(M_n)}
\def \Oezm {\Cal{O}_\varepsilon^{\,\Zeps\!}(M_n)}

\def \Ounozm {\Cal{O}_1^{\,\Bbb{Z}}(M_n)}
\def \Om {\Cal{O}(M_n)}
\def \Ozm {\Cal{O}^{\,\Bbb{Z}}(M_n)}
\def \Ozem {\Cal{O}^{\,\Zeps\!}(M_n)}

\def \B {\Bbb{B}}
\def \Bb {\text{\rm B}}
\def \t {\underline{t}}
\def \gerFr {{\frak F}{\frak r}}

\document


\topmatter

{\ } 

\vskip-33pt  

\hfill   
  {\smallrm {\smallsl Glasgow Mathematical Journal\/} 
{\smallbf 49}  (2007), 479--488}   
\hskip19pt   {\ }  

\vskip51pt  

\title
  PBW theorems and Frobenius structures for quantum matrices  
\endtitle

\author
       Fabio Gavarini
\endauthor

\leftheadtext{ Fabio Gavarini }
\rightheadtext{ PBW theorems and Frobenius structures for
quantum matrices }

\affil
  Universit\`a di Roma ``Tor Vergata'' ---
Dipartimento di Matematica  \\
  Via della Ricerca Scientifica 1, I-00133 Roma --- ITALY  
\endaffil

\address\hskip-\parindent
  Fabio Gavarini  \newline   
  \indent   Universit\`a degli Studi di Roma ``Tor Vergata''  
---   Dipartimento di Matematica  \newline   
  \indent   via della ricerca scientifica 1,
I-00133 Roma, ITALY  
---   gavarini\@{}mat.uniroma2.it  \newline   
  \indent   {\tt http://www.mat.uniroma2.it/\~{}gavarini/page-web.html}   
\endaddress   

\abstract
   Let  $ \, G \in \big\{ \text{\it Mat}_n(\C) , {GL}_n(\C) , {SL}_n(\C)
\big\} \, $,  \, let  $ \Oqg $  be the quantum function algebra   ---
over  $ \Z\big[q,q^{-1}\big] $  ---   associated to  $ G $,  and let 
$ \Oeg $  be the specialisation of the latter at a root of unity 
$ \varepsilon $,  whose order  $ \ell $  is odd.  There is a quantum
Frobenius morphism that embeds  $ \Og $,  the function algebra of 
$ G $,  in  $ \Oeg $  as a central Hopf subalgebra, so that  $ \Oeg $ 
is a module over  $ \Og $.  When  $ \, G = {SL}_n(\C) \, $,  \, it is
known by [BG], [BGStaf] that (the complexification of) such a module
is free, with rank  $ \ell^{\text{\it dim}(G)} \, $.  In this note we
prove a PBW-like theorem for  $ \Oqg $,  \, and we show that   --- when 
$ G $  is  $ \text{\it Mat}_{\,n} $  or  $ {GL}_n $  ---   it yields
explicit bases of  $ \Oeg $  over  $ \Og \, $.  As a direct application,
we prove that  $ \Oegl $  and  $ \Oem $  are free Frobenius extensions
over  $ \Ogl $  and  $ \Om \, $,  thus extending some results of
[BGStro].   
\endabstract

\endtopmatter

\footnote""{Keywords: \ {\sl PBW theorems, quantum groups,
roots of unity}.}

\footnote""{2000 {\it Mathematics Subject Classification:}
\ Primary 20G42; Secondary 81R50.}

\vskip9pt

\centerline {\bf \S \; 1 \, The general setup }

\vskip13pt

   Let  $ G $  be a complex semisimple, connected, simply connected
affine algebraic group.  One can introduce a quantum function
algebra  $ \Oqg $,  a Hopf algebra over the ground ring  $ \C \!\left[
q, q^{-1}\right] $,  where  $ q $  is an indeterminate,  as in [DL].  If 
$ \varepsilon $  is any root of 1, one can specialize  $ \Oqg $  at  $ \,
q = \varepsilon \, $,  which means taking the Hopf  $ \C $--algebra  $ \;
\Oeg := \Oqg \Big/ (q-\varepsilon) \Oqg \; $.  In particular, for  $ \,
\varepsilon = 1 \, $  one has  $ \; \Ounog \cong \Og \, $,  \; the
classical (commutative) function algebra over  $ G \, $.  Moreover, if
the order  $ \ell $  of  $ \varepsilon $  is odd, then there exists a
Hopf algebra monomorphism  $ \; \gerFr \, \colon \,  \Og \cong \Ounog
\!\llonghookrightarrow \! \Oeg \, $,  \; called  {\it quantum Frobenius
morphism for\/}  $ G \, $,  which embeds  $ \Og $  inside  $ \Oeg $  as
a central Hopf subalgebra.  Therefore,  $ \Oeg $  is naturally a module
over  $ \Og \, $.  It is proved in [BGStaf] and in [BG] that such a
module is free, with rank  $ \ell^{\text{\it dim}(G)} \, $.  In the
special case of  $ \, G = {SL}_2 \, $,  a stronger result was given in
[DRZ], where an explicit basis was found.  We shall give similar results
when  $ G $  is  $ {GL}_n $  or  $ \, M_n := \text{\it Mat}_n \, $;  \,
namely we provide explicit bases of  $ \Oeg $  as a free module over 
$ \Og \, $,  where in addition everything is defined replacing  $ \C $ 
with  $ \Z \, $.  The proof is via some (more or less known) PBW
theorems for  $ \Oqm \, $  and  $ \Oqgl $   --- and  $ \Oqsl $ 
as well ---   as modules over  $ \Z \big[ q, q^{-1} \big] \, $.  

\vskip7pt

   Let  $ \, M_n := \text{\it Mat}_n(\C) \, $.  The algebra  $ \Om $ 
of regular functions on  $ M_n $  is the unital associative commutative 
$ \C $--algebra  with generators  $ \, \bar{t}_{i,j} \, $  ($ \, i, j = 1,
\dots, n \, $).  The semigroup structure on  $ M_n $  yields on  $ \Om $ 
the natural bialgebra structure given by matrix product   --- see [CP],
Ch.~7.  We can also consider the semigroup-scheme  $ \big( M_n \big)_\Z $ 
associated to  $ M_n $,  \, for which a like analysis applies: in
particular, its function algebra  $ \Ozm $  is a  $ \Z $--bialgebra, 
with the same presentation as  $ \Om $  but over the ring  $ \Z \, $.   

\vskip7pt

   Now we define quantum function algebras.  Let  $ R $  be any
commutative ring with unity, and let  $ \, q \in R \, $  be invertible. 
We define  $ \Oqrm $  as the unital associative  $ R $--algebra 
with generators  $ \; t_{i,j} \; $  ($ \, i, j = 1, \dots, n \, $) 
and relations   
 \vskip-13pt   
  $$  \hbox{ $ \eqalign {
   {} \hskip11pt   t_{i,j} \, t_{i,k} = q \, t_{i,k} \, t_{i,j} \; ,
\quad \quad  t_{i,k} \, t_{h,k}  &  = q \, t_{h,k} \, t_{i,k}
\hskip56pt  \forall  \quad  j<k \, , \, i<h \, ,  \cr
   {} \hskip11pt   t_{i,l} \, t_{j,k} = t_{j,k} \, t_{i,l} \; ,
\qquad  t_{i,k} \, t_{j,l} - \, t_{j,l} \, t_{i,k}  &
= \left( q - q^{-1} \right) \, t_{i,l} \, t_{j,k}   \hskip20pt
\forall  \quad  i<j \, , \, k<l \, .  \cr } $ }  $$
   \indent   It is known that  $ \Oqrm $  is a bialgebra, but we do
not need this extra structure in the present work (see [CP] for further
details   ---   cf.~also [AKP] and [PW]).   

\vskip5pt   

   As to specialisations, set  $ \, \Z_q := \Z \big[ q, q^{-1}\big] \, $, 
\, let  $ \, \ell \in \N_+ \, $  be odd, let  $ \, \phi_\ell(q) \, $  be
the  $ \ell $-th  cyclotomic polynomial in  $ q \, $,  \, and let 
$ \, \varepsilon := \overline{q} \in \Zeps := \Z_q \Big/ \big(
\phi_\ell(q) \big) \, $,  \, so that  $ \varepsilon $  is a (formal)
primitive  $ \ell $-th  root of 1 in  $ \Zeps \, $.  Then   
 \vskip-5pt
  $$  \Oezm  \; = \;  \Oqzm \Big/ \big(\phi_\ell(q)\big) \, \Oqzm 
\; \cong \;  \Zeps \otimes_\Z \Oqzm  \quad .  $$   
 \vskip-1pt
   It is also known that there is a bialgebra isomorphism   
 \vskip-13pt
  $$  \Ounozm \cong \Oqzm \! \Big/ \! (q \! - \! 1) \Oqzm \,
\lhook\joinrel\relbar\joinrel\twoheadrightarrow \Ozm \; ,  \quad
\! t_{i,j} \!\!\! \mod (q \! - \! 1) \, \Oqzm \; \mapsto \;
\bar{t}_{i,j}  $$   
 \vskip-6pt
\noindent   
 and a bialgebra monomorphism, called  {\it quantum
Frobenius morphism\/}  ($ \varepsilon $  and  $ \ell \, $  as above),  
  $$  \gerFr_{\Z} \; \colon \; \Ozm \cong \Ounozm \,
\llonghookrightarrow \, \Oezm \;\; ,  \qquad  \bar{t}_{i,j}
\mapsto t_{i,j}^{\,\ell}{\big|}_{q=\varepsilon}  $$   
whose image is  {\sl central\/}  in  $ \Oezm \, $.  Thus  $ \, \Ozem :=
\Zeps \otimes_\Z \Ozm \, $  becomes identified   --- via  $ \gerFr_{\Z} \, $, 
which clearly extends to  $ \Ozem $  by scalar extension ---   with a
central subbialgebra of  $ \Oezm \, $,  \, so the latter can be seen
as an  $ \Ozem $--module.  By the result in [BGStaf] and [BG] mentioned
above, we can expect this module to be free, with rank  $ \ell^{n^2}
\, $.   

\vskip5pt   

   All the previous framework also extends to  $ \, {GL}_n \, $  and
to  $ \, {SL}_n \, $  instead of  $ \, M_n \, $.  Indeed, consider
the  {\it quantum determinant}  $ \; D_q := \sum_{\sigma \in
\Cal{S}_n} {(-q)}^{\ell(\sigma)} t_{1,\sigma(1)} \, t_{2,\sigma(2)}
\cdots t_{n,\sigma(n)} \in \Oqrm \, $,  \; where  $ \ell(\sigma) $ 
denotes the  {\sl length\/}  of any permutation  $ \sigma $  in the
symmetric group  $ \Cal{S}_n \, $.  Then  $ D_q $  belongs to the
centre of  $ \Oqrm $,  hence one can extend  $ \Oqrm $  by a formal
inverse to  $ D_q \, $,  \, i.e.~defining the algebra  $ \, \Oqrgl
:= \Oqrm\big[D_q^{-1}\big] \, $.  Similarly, we can define also  $ \,
\Oqrsl := \Oqrm \Big/ \big( D_q - 1 \big) \, $.  Now  $ \Oqrgl $  and 
$ \Oqrsl $  are  {\sl Hopf}  $ R $--algebras,  and the maps  $ \, \Oqrm
\longhookrightarrow \Oqrgl \, $,  $ \, \Oqrgl \longtwoheadrightarrow
\Oqrsl \, $,  $ \, \Oqrm \longtwoheadrightarrow \Oqrsl \, $  (the
third one being the composition of the first two) given by  $ \, t_{i,j}
\mapsto t_{i,j} \, $  are epimorphisms of  $ R $--bialgebras,  and even
of Hopf  $ R $--algebras  in the second case.  The specialisations   
 \vskip-13pt
  $$  \displaylines{ 
   \Oezgl  \; = \;  \Oqzgl \Big/ \big(\phi_\ell(q)\big) \, \Oqzgl 
\; \cong \;  \Zeps \otimes_\Z \Oqzgl  \cr   
   \Oezsl  \; = \;  \Oqzsl \Big/ \big(\phi_\ell(q)\big) \, \Oqzsl 
\; \cong \;  \Zeps \otimes_\Z \Oqzsl  \cr }  $$   
 \vskip-5pt
\noindent   
 enjoy the same properties as above, namely there exist isomorphisms 
$ \; \Ounozgl \cong \Ozgl \; $  and  $ \; \Ounozsl \cong \Ozsl \; $ 
and there are quantum Frobenius morphisms 
 \vskip-15pt
  $$  \gerFr_{\Z} \; \colon \, \Ozgl \cong \Ounozgl \, \longhookrightarrow
\, \Oezgl \; ,  \quad  \gerFr_{\Z} \; \colon \, \Ozsl \cong \Ounozsl
\, \longhookrightarrow \, \Oezsl  $$   
%
%
 \eject   
\noindent   
described by the same formul{\ae}  as for  $ M_n \, $.  Moreover,
$ \; D_q^{\pm 1} \!\! \mod (q-1) \, \mapsto \, D^{\pm 1} \; $  in the
isomorphisms and  $ \; D^{\pm 1} \cong D_q^{\pm 1} \!\! \mod (q-1) \,
\mapsto \, D_q^{\pm \ell} \!\! \mod (q-\varepsilon) \; $  in the quantum
Frobenius morphisms for  $ {GL}_n $  (which extend those of  $ M_n \, $). 
In addition, all these isomorphisms and quantum Frobenius morphisms are
compatible (in the obvious sense) with the natural maps which link 
$ \Oqzm $,  $ \Oqzgl $  and  $ \Oqzsl $,  and their specialisations,
to each other.  
                                                 \par   
   Like for  $ M_n \, $,  \, the image of the quantum Frobenius
morphisms  are  {\sl central\/}  in  $ \Oezgl $  and in  $ \Oezsl \, $. 
Thus  $ \, \Ozegl := \Zeps \otimes_\Z \Ozgl \, $  identifies to a central
Hopf subalgebra of  $ \Oezgl \, $,  and  $ \, \Ozesl := \Zeps \otimes_\Z
\Ozsl \, $  identifies to a central Hopf subalgebra of  $ \Oezsl \, $; 
\, so  $ \Oezgl $  is an  $ \Ozgl $--module  and  $ \Oezsl $  is an 
$ \Ozsl $--module.   

\vskip1pt   

   In \S 2, we shall prove (Theorem 2.1) a PBW-like theorem providing
several different bases for  $ \Oqrm $,  $ \Oqrgl $  and  $ \Oqrsl $ 
as  $ R $--modules.  As an application, we find (Theorem 2.2) explicit
bases of  $ \Oezm $  as an  $ \Ozem $--module,  which then in particular
is free of rank  $ \ell^{\text{\it dim}(M_n)} \, $.  The same bases are
also  $ \Ozegl $--bases  for  $ \Oezgl $,  which then is free of rank 
$ \ell^{\text{\it dim}({GL}_n)} \, $.  Both results can be seen as
extensions of some results in [BGStaf].   

\vskip1pt   

   Finally, in \S 3 we use the above mentioned bases to prove that 
$ \Ozem $  is a  {\sl free Frobenius extension\/}  of its central
subalgebra  $ \Ozem $,  and to explicitly compute the associated
Nakayama automorphism.  The same we do for  $ \Oezgl $  as well. 
Everything follows from the ideas and methods in [BGStro], now
applied to the explicit bases given by Theorem 2.2.   

\vskip27pt

\centerline { \bf  \S \; 2 \, PBW--like theorems }

\vskip13pt

\proclaim{Theorem 2.1}  {\sl (PBW theorem for  $ \Oqrm $,  $ \Oqrgl $ 
and  $ \Oqrsl $  as  $ R $--modules)}   
                                                    \par   
\noindent
   \hskip3pt   Assume  $ (q-1) $  {\sl is not invertible in} 
$ \, R_q := \big\langle q, q^{-1} \big\rangle \, $,  \, the subring
of  $ \, R $  generated by  $ q $  and  $ q^{-1} $.     
 \vskip3pt   
   (a) \, Let any total order be fixed in  $ \, {\{ 1, \dots, n
\}}^{\times 2} \, $.  Then the following sets of ordered monomials
are  $ R $--bases  of  $ \, \Oqrm \, $,  \, resp.~$ \, \Oqrgl \, $, 
\, resp.~$ \, \Oqrsl \, $,  \, as modules over  $ R \, $:   
 \vskip-9pt   
  $$  \displaylines{ 
   B_M \;  :=  \; \bigg\{\, {\textstyle \prod\limits_{i,j=1}^n}
t_{i,j}^{N_{i,j}} \;\bigg\vert\; N_{i,j} \in \N \; \forall \,
i, j \,\bigg\}  \cr    
   B_{GL}^{\,\wedge} \;  :=  \; \bigg\{\, {\textstyle
\prod\limits_{i,j=1}^n} t_{i,j}^{N_{i,j}} \, D_q^{-N}
\;\bigg\vert\; N, N_{i,j} \in \N \; \forall \, i, j \; ;
\,\; \min \big( \big\{ N_{i,i} \big\}_{1 \leq i \leq n} \!
\cup \{N\} \big) = 0 \,\bigg\}  \cr   
   B_{GL}^{\,\vee} \;  :=  \; \bigg\{\, {\textstyle \prod\limits_{i,j=1}^n}
t_{i,j}^{N_{i,j}} \, D_q^Z \;\bigg\vert\; Z \in \Z \, , N_{i,j} \in
\N \; \forall \, i, j \; ; \,\; \min \big\{ N_{i,i} \big\}_{1 \leq i
\leq n} = 0 \,\bigg\}  \cr   
   B_{SL} \;  :=  \; \bigg\{\, {\textstyle \prod\limits_{i,j=1}^n}
t_{i,j}^{N_{i,j}} \;\bigg\vert\; N_{i,j} \in \N \; \forall \, i, j \; ;
\,\; \min \big\{ N_{i,i} \big\}_{1 \leq i \leq n} = 0 \,\bigg\} }  $$   
 \vskip3pt   
   (b) \, Let  $ \preceq $  be any total order fixed in  $ \, {\{ 1,
\dots, n \}}^{\times 2} $  such that  $ \, (i,j) \preceq (h,k) \preceq
(l,m) \, $  whenever  $ \, j > n \! + \! 1 \! - \! i \, $,  $ \, k =
n \! + \! 1 \! - \! h \, $,  $ \, m < n \! + \! 1 \! - \! l \, $.  Then
the following sets of ordered monomials are  $ R $--bases  of  $ \,
\Oqrgl \, $,  \, resp.~$ \, \Oqrsl \, $,  \, as modules over  $ R \, $:   
 \vskip-9pt   
  $$  \displaylines{ 
   B^{\wedge,-}_{GL} \;  :=  \; \bigg\{\, {\textstyle \prod\limits_{i,j=1}^n}
t_{i,j}^{N_{i,j}} \, D_q^{-N} \;\bigg\vert\; N, N_{i,j} \in \N \; \forall
\, i, j \; ; \,\; \min \big( \big\{ N_{i,n+1-i} \big\}_{1 \leq i \leq n} \!
\cup \{N\} \big) = 0 \,\bigg\}  \cr   
   B^{\vee,-}_{GL} \;  :=  \; \bigg\{\, {\textstyle \prod\limits_{i,j=1}^n}
t_{i,j}^{N_{i,j}} \, D_q^Z \;\bigg\vert\; Z \in \Z \, , N_{i,j} \in \N
\; \forall \, i, j \; ; \,\; \min \big\{ N_{i,n+1-i} \big\}_{1 \leq i
\leq n} = 0 \,\bigg\}  \cr   
   B^-_{SL} \;  :=  \; \bigg\{\, {\textstyle \prod\limits_{i,j=1}^n}
t_{i,j}^{N_{i,j}} \;\bigg\vert\; N_{i,j} \in \N \; \forall \, i, j \; ;
\,\; \min \big\{ N_{i,n+1-i} \big\}_{1 \leq i \leq n} = 0 \,\bigg\} }  $$  
\endproclaim   
 \eject   

\demo{Proof}  Roughly speaking, our method is a (partial) application
of the diamond lemma (see [Be]): however, we do not follow it in
all details, as we use a specialisation trick as a shortcut.   
                                                 \par   
   If we prove our results for the algebras defined over  $ R_q $ 
instead of  $ R \, $,  then the same results will hold as well by
scalar extension.  Thus we can assume  $ \, R = R_q \, $,  \, and
then we note that, by our assumption,  {\sl the specialised ring  $ \,
\overline{R} := R \big/ (q-1) R \not= \{0\} \, $  is non-trivial}.  

\vskip9pt

   {\it  $ \underline{\text{Proof of (a)}} $:} \,  {\sl (see also
[Ko], Theorem 3.1, and [PW], Theorem 3.5.1)}   
 \vskip3pt
   We begin with  $ \Oqrm \, $.  It is clearly spanned over  $ R $ 
by the set of all (possibly unordered) monomials in the  $ t_{ij} $'s: 
so we must only prove that any such monomial belongs to the  $ R $--span 
of the  {\sl ordered\/}  monomials.  In fact, the latter are linearly
independent, since such are their images via specialisation  $ \; \Oqrm
\longtwoheadrightarrow \Oqrm \Big/ (q \! - \! 1) \, \Oqrm \, \cong \,
\Cal{O}^{\,\overline{R}}_1(M_n) \; $.   
                                             \par    
   Thus, take any (possibly unordered) monomial in the  $ t_{ij} $'s, 
say  $ \; \t \, := \, t_{i_1,j_i} \, t_{i_2,j_2} \cdots \,
t_{i_k,j_k} \; $,  \; where  $ k $  is the degree of  $ \t
\, $:  \, we associate to it its  {\sl weight},  defined as     
  $$  w(\,\t\,) \; := \; \big( k \, , d_{1,1} \, , d_{1,2} \, ,
\dots, d_{1,n} \, , d_{2,1} \, , d_{2,2} \, , \dots, d_{2,n} \, , d_{3,1}
\, , \dots, d_{n-1,n} \, , d_{n,1} \, , d_{n,2} \, , \dots, d_{n,n}
\big)  $$   
where  $ \; d_{i,j} := \big| \big\{\, s \! \in \! \{1,\dots,k\} \,\big|\,
(i_s,j_s) = (i,j) \big\} \big| \, $  = {\sl number of occurrences of 
$ t_{i,j} $  in  $ \t \, $}.  Then  $ \, w(\t) \in \N^{n^2 + 1} \, $, 
\, and we consider  $ \N^{n^2 + 1} $  as a totally ordered set with
respect to the (total) lexicographic order  $ \leq_{lex} \, $.  By
a quick look at the defining relations of  $ \Oqrm $,  namely    
  $$  \eqalignno {
   {} \hskip11pt   t_{i,j} \, t_{i,k} = q \, t_{i,k} \, t_{i,j} \; ,
\quad \quad  t_{i,k} \, t_{h,k}  &  = q \, t_{h,k} \, t_{i,k}  \hskip51pt 
&  \forall  \quad  j<k \, , \, i<h \, ,  \cr
   {} \hskip11pt   t_{i,l} \, t_{j,k} = t_{j,k} \, t_{i,l} \; ,
\qquad  t_{i,k} \, t_{j,l} - \, t_{j,l} \, t_{i,k}  &
= \left( q - q^{-1} \right) \, t_{i,l} \, t_{j,k}  \hskip51pt  &  
\forall  \quad  i<j \, , \, k<l \, .  \cr }  $$
one easily sees that the weight defines an  {\sl algebra filtration\/} 
on  $ \Oqrm \, $.
                                           \par   
   Now, using these same relations, one can re-order the  $ t_{ij} $'s 
in any monomial according to the fixed total order.  During this process,
only two non-trivial things may occur, namely:  
 \vskip5pt   
  {\sl --1)} \, some powers of  $ q $  show up as coefficients
(when a relation in first line is employed);  
                                             \par    
  {\sl --2)} \, a new summand is added (when the bottom-right
relation is used);  
 \vskip5pt   
   If only steps of type  {\sl 1)\/}  occur, then the process eventually
stops with an  {\sl ordered\/}  monomial in the  $ t_{ij} $'s  multiplied
by a power of  $ q \, $.  Whenever instead a step of type  {\sl 2)\/} 
occurs, the newly added term is just a coefficient  $ \big( q - q^{-1}
\big) $  times a (possibly unordered) monomial in the  $ t_{ij} $'s, 
call it  $ \t' \, $:  \, however, by construction  $ \, w
\big( \t'\big) \lneqq_{lex} w(\,\t\,) \, $. 
Then, by induction on the weight, we can assume that  $ \t' $ 
lies in the  $ R $--span  of the ordered monomials, so we can ignore
the new summand.  The process stops in finitely many steps, and we are
done with  $ \Oqrm \, $.
 \vskip5pt
   Second, we look at  $ \Oqrgl \, $.  Let us consider  $ \, f \in \Oqrgl
\, $.  By definition, there exists  $ \, N \in \N \, $  such that  $ \,
f D_q^N \in \Oqrm \, $;  \, therefore, by the result for  $ \Oqrm $ 
just proved, we can expand  $ \, f D_q^N \, $  as an  $ R $--linear 
combination of ordered monomials, call them  $ \, \t =
\prod_{i,j=1}^n t_{i,j}^{N_{i,j}} \, $.  Thus,  $ f $  itself
is an  $ R $--linear  combination of monomials  $ \; \t
\, D_q^{-N} \, $,  \, so the latter span  $ \Oqrgl \, $.
                                          \par   
   Now consider an ordered monomial  $ \; \t = \prod_{i,j=1}^n
t_{i,j}^{N_{i,j}} \, $  in which  $ \, N_{i,i} > 0 \, $  for all 
$ i \, $.  Then we can re-arrange the  $ t_{i,i} $'s  in  $ \t \, $ 
so to single out a factor  $ \; t_{1,1} \, t_{2,2} \, \cdots \,
t_{n-1,n-1} \, t_{n,n} \, $,  \; up to ``paying the cost'' (perhaps)
of producing some new summands  {\sl of lower weight\/}:  the outcome
reads
  $$  \t \, = \; q^s \, \t_{\,0} \, t_{1,1} \, t_{2,2}
\, \cdots \, t_{n-1,n-1} \, t_{n,n} \, + \; \text{\it l.t.'s}  
\eqno (2.1)  $$   
for some  $ \, s \in \Z \, $,  \, with  $ \, \t_{\,0} := \prod_{i,j=1}^n
t_{i,j}^{N_{i,j} - \delta_{i,j}} \, $  having lower weight than  $ \t
\, $,  \, and the expression  {\it l.t.'s\/}  standing for an 
$ R $--linear  combination of some monomials  $ \; \underline{\check{t}}
\; $  such that  $ \, w \big( \underline{\check{t}} \,\big) \lneqq_{lex}
w \big( \t \big) \, $.  Then we re-write the monomial  $ \, t_{1,1} \,
t_{2,2} \, \cdots \, t_{n-1,n-1} \, t_{n,n} \, $  using the identity   
  $$  t_{1,1} \, t_{2,2} \, \cdots \, t_{n-1,n-1} \, t_{n,n}  \; = \; 
D_q  -  {\textstyle \sum_{\Sb  \sigma \in \Cal{S}_n  \\   \sigma \not=
\text{\it id}  \endSb}}  \!\! {(-q)}^{\ell(\sigma)} \, t_{1,\sigma(1)}
\, t_{2,\sigma(2)} \cdots t_{n,\sigma(n)}  \; = \; 
D_q  +  \text{\it l.t.'s}   \eqno (2.2)  $$   
 \vskip-5pt   
\noindent   
 and we replace the right-hand side of (2.2) inside (2.1).  We get  $ \; \t
\, = \, q^s \, \t_{\,0} \, D_q + \, \text{\it l.t.'s} \; $  (for  $ D_q $ 
is central!), where now  $ \; \t_{\,0} \; $  and all monomials within 
{\it l.t.'s\/}  have strictly lower weight than  $ \, \t \, $.  
                                                   \par   

 If we look now at  $ \, \t \, D_q^{\,z} \, $  (for some  $ \, z \in \Z
\, $),  we can re-write  $ \t $  as above, thus getting   
  $$  \t \, D_q^z  \; = \;  q^s \, \t_{\,0} \, D_q \, D_q^{\,z}
\, + \;  \text{\it l.t.'s}  \; = \;  q^s \, \t_{\,0} \, D_q^{\,z+1}
\, + \;  \text{\it l.t.'s}   \eqno   (2.3)  $$   
where  {\it l.t.'s\/}  is an  $ R $--linear  combination of monomials 
$ \; \underline{\tilde{t}} \, D_q^{\,z+1} \; $  such that  $ \, w \big(
\underline{\tilde{t}} \,\big) \lneqq_{lex} w \big( \t \big)
\; $.
                                                   \par   
   By repeated use of (2.3) as reduction argument, we can easily
show   --- by induction on the weight ---   that any monomial of
type  $ \, \t \, D_q^{-N} \, $  ($ N \in \N \, $)  can be expanded as
an  $ R $--linear  combination elements of  $ B^{\,\wedge}_{GL} $  or
elements of  $ B^{\,\vee}_{GL} \, $.  Thus, both these sets do span 
$ \Oqrgl \, $.   
                                                   \par   
   To finish with, both  $ B^{\,\wedge}_{GL} $  and 
$ B^{\,\vee}_{GL} $  are  $ R $--linearly  independent,
as their image through the specialisation epimorphism 
$ \, \Oqrgl \! \relbar\joinrel\twoheadrightarrow \!
\Cal{O}_1^{\,\overline{R}}({GL}_n) \cong \Cal{O}^{\,\overline{R}}({GL}_n)
\, $  are  $ \overline{R} $--bases  of  $ \Cal{O}^{\,\overline{R}}({GL}_n)
\, $.   
 \vskip5pt
   As to  $ \Oqrsl \, $,  we can repeat the argument for  $ \Oqrgl \, $. 
First,  $ B_{SL} $  is linearly independent, for its image through
specialisation  $ \; \Oqrsl \relbar\joinrel\twoheadrightarrow
\Cal{O}_1^{\,\overline{R}}({SL}_n) \cong \Cal{O}^{\,\overline{R}}
({SL}_n) \; $  is an  $ \overline{R} $--basis  of  $ \Cal{O}^{\,
\overline{R}}({SL}_n) \, $.  Second, the epimorphism  $ \; \Oqrm
\relbar\joinrel\twoheadrightarrow \Oqrsl \; \big( t_{i,j} \mapsto
t_{i,j} \big) \, $,  and the result for  $ \Oqrm \, $,  \, imply
that the  $ R $--span  of  $ \, S_{SL} := \bigg\{\, {\textstyle
\prod\limits_{i,j=1}^n} t_{i,j}^{N_{i,j}} \;\bigg\vert\; N_{i,j}
\! \in \! \N \; \forall \, i, j \,\bigg\} \, $  is  $ \Oqrsl \, $. 
Thus one is only left to prove that each monomial  $ \, \t =
\prod_{i,j=1}^n t_{i,j}^{N_{i,j}} \in S_{SL} \, $  belongs to
the  $ R $--span  of  $ B_{SL} \, $:  as before, this can be done
by induction on the weight, using the reduction formula  $ \; \t
\, = \, q^s \, \t_{\,0} \, D_q + \, \text{\it l.t.'s} \; $  (see
above), and plugging in it the relation  $ \, D_q = 1
\, $.   
                                               \par   
   Alternatively, we remind there is an isomorphism 
$ \, \Oqrsl \otimes_R R\big[x,x^{-1}\big] \cong \Oqrgl \, $ 
(of  $ R $--algebras)  given by  $ \; t_{i,j} \otimes x^z \mapsto
D_q^{-\delta_{i,1}} \, t_{i,j} \cdot D_q^{\,z} \, $  (cf.~[LS]). 
This along with the result about  $ B^\vee_{GL} $  clearly implies
that also  $ B_{SL} $  is an  $ R $--basis  for  $ \Oqrsl \, $, 
\, as claimed.   

\vskip9pt

   {\it  $ \underline{\text{Proof of (b)}} $:} \,  First look at 
$ \Oqrgl \, $.  If  $ \, f \in \Oqrgl \, $,  like in the proof of 
{\it (a)\/}  we expand  $ \, f D_q^N \, $  as an  $ R $--linear 
combination of ordered (according to  $ \preceq \, $)  monomials of
type  $ \; \t =  \t^- \, \t^= \, \t^+ \, $,  \, with  $ \; \t^- :=
\prod_{j>n+1-i} t_{i,j}^{N_{i,j}} \, $,  $ \; \t^= := \prod_{j=n+1-i}
t_{i,j}^{N_{i,j}} \, $  and  $ \; \t^+ := \prod_{j<n+1-i}
t_{i,j}^{N_{i,j}} \; $.  So  $ f $  is an  $ R $--linear 
combination of monomials  $ \; \t^- \, \t^= \, \t^+ \,
D_q^{-N} \, $,  \, hence the latter span  $ \Oqrgl \, $.   
                                          \par   
   We show that each (ordered) monomial  $ \; \t^- \, \t^= \, \t^+ \,
D_q^{-N} \; $  belongs both to the  $ R $--span  of  $ B^{\wedge,-}_{GL} $ 
and of  $ B^{\vee,-}_{GL} \, $,  \, by induction on the (total) degree of
the monomial  $ \, \t^= \, $.  The basis of induction is  $ \; \text{\it
deg}\,(\t^=) = 0 \, $,  \, so that  $ \, \t^= = 1 \, $  and  $ \; \t^-
\, \t^= \, \t^+ \, D_q^{-N} = \, \t^- \, \t^+ \, D_q^{-N} \in \,
B^{\wedge,-}_{GL} \cap B^{\vee,-}_{GL} \; $.   
                                          \par   
   As a matter of notation, let  $ \Cal{N}^- $,  {resp.}  $ \Cal{H} \, $, 
resp.~$ \Cal{N}^+ $,  be the  $ R $--subalgebra  of  $ \Oqrm $  generated
by the  $ t_{i,j} $'s  with  $ \, j > n \! + \! 1 \! - i \, $,  {resp.}
$ \, j = n \! + \! 1 \! - i \, $,  resp.~$ \, j < n \! + \! 1 \! - i \, $. 
Note that  $ \Cal{H} $  is Abelian, and  $ \, \t^- \in \Cal{N}^- \, $, 
$ \, \t^= \in \Cal{H} \, $,  $ \, \t^+ \in \Cal{N}^+ \, $.   
                                          \par   
   Now assume that all the exponents  $ N_{i,n+1-i} $'s  in the factor 
$ \t^= $  are strictly positive.  As  $ \Cal{H} $  is Abelian, we can
draw out of  $ \t^= $  (even out of  $ \, \t =  \t^- \, \t^= \, \t^+
\, $)  a factor  $ \; t_{n,1} \, t_{n-1,2} \, \cdots \, t_{2,n-1} \,
t_{1,n} \; $.  Now recall that  $ D_q $  can be expanded as  $ \,
D_q = {\textstyle \sum_{\sigma \in \Cal{S}_n}} \! {(\!-q)}^{\ell(\sigma)}
t_{n,\sigma(n)} \, t_{n-1,\sigma(n-1)} \cdots t_{2,\sigma(2)} \,
t_{1,\sigma(1)} $  \, (see, e.g., [PW] or [Ko]).  Then we can
re-write the monomial  $ \; t_{n,1} \, t_{n-1,2} \, \cdots \,
t_{2,n-1} \, t_{1,n} \; $  as   
  $$  t_{n,1} \, t_{n-1,2} \, \cdots \, t_{1,n}  \; = \, 
{(-q)}^{-\ell(\sigma_0)} D_q  -  {\textstyle \sum_{\Sb  \sigma
\in \Cal{S}_n  \\   \sigma \not= \sigma_0  \endSb}}  \!\!
{(-q)}^{\ell(\sigma) - \ell(\sigma_0)} \, t_{n,\sigma(n)} \,
t_{n-1,\sigma(n-1)} \cdots t_{1,\sigma(1)}   \eqno (2.4)  $$   
where  $ \, \sigma_0 \in \Cal{S}_n \, $  is the permutation  $ \,
i \mapsto (n+1-i) \, $.  Note also that we can reorder the factors
in the summands of (2.4) so that all factors  $ t_{i,j} $  from 
$ \Cal{N}^- $  are on the left of those from  $ \Cal{N}^+ $. 
                                                \par   
   Now we replace the right-hand side of (2.4) in the factor 
$ \t^= \, $  within  $ \, \t =  \t^- \, \t^= \, \t^+ \, $, 
\, thus   
 \vskip4pt   
   \centerline{ $ \t^- \, \t^= \, \t^+  \; = \; 
{(-q)}^{-\ell(\sigma_0)} \, \t^- \, \t^=_{\,0} \,
D_q \, \t^+ \, + \, \text{\it l.t.'s}  \; = \; 
{(-q)}^{-\ell(\sigma_0)} \, \t^- \, \t^=_{\,0} \,
\t^+ \, D_q \, + \, \text{\it l.t.'s} $ }   
 \vskip3pt   
\noindent   
 Here  $ \, \t^=_{\,0} := \t^= \, {\big( t_{n,1} \, t_{n-1,2} \,
\cdots \, t_{2,n-1} \, t_{1,n} \big)}^{-1} \, $  has lower (total)
degree than  $ \t^= \, $,  \, and the expression  {\it l.t.'s\/} 
stands for an  $ R $--linear  combination of some other monomials 
$ \; \underline{\hat{t}}^{\,-} \, \underline{\hat{t}}^{\,=} \,
\underline{\hat{t}}^{\,+} \; $  (like  $ \; \t^- \, \t^= \, \t^+ \; $ 
above) in which again the degree of  $ \, \underline{\hat{t}}^{\,=} $ 
is lower than the degree of  $ \, \t^= \, $.  In fact, this holds because
when any factor  $ \, t_{i,\sigma(i)} \in \Cal{N}^- \, $  is pulled from
the right to the left of any monomial in  $ \, \underline{\check{t}}^{\,=}
\in \Cal{H} \, $  the degree of  $ \underline{\check{t}}^{\,=} $  is not
increased.  By induction on this degree, we can easily conclude that
every ordered monomial  $ \, \t^- \, \t^= \, \t^+ D_q^{\,z} \, $ 
(with  $ \, z \in \Z \, $)  belongs to both the  $ R $--span  of 
$ B^{\wedge,-}_{GL} $  and the  $ R $--span  of  $ B^{\vee,-}_{GL} \, $. 
That is, both sets span  $ \Oqrgl \, $.  
                                                   \par   
   Eventually, both  $ B^{\wedge,-}_{GL} $  and  $ B^{\vee,-}_{GL} $ 
are linearly independent, as their image through the specialisation
epimorphism  $ \; \Oqrgl \longtwoheadrightarrow \Cal{O}_1^{\,\overline{R}}
({GL}_n) \cong \Cal{O}^{\,\overline{R}}({GL}_n) \; $  are 
$ \overline{R} $--bases  of  $ \Cal{O}^{\,\overline{R}}({GL}_n) \, $.   
 \vskip5pt
   Second, we look at  $ \Oqrsl \, $.  Like for claim  {\it (a)},  we
can repeat again   --- {\it mutatis mutandis}  ---   the argument for 
$ \Oqrgl \, $,  \, which does work again   ---   one only has to plug
in the additional relation  $ \, D_q = 1 \, $  too.  Otherwise, as
an alternative proof, we can note that the isomorphism  $ \, \Oqrsl
\otimes_R R\big[x,x^{-1}\big] \cong \Oqrgl \, $  together with the
result about  $ B^{\vee,-}_{GL} $  easily implies that  $ B^-_{SL} $ 
too is an  $ R $--basis  for  $ \Oqrsl \, $,  \, q.e.d.   \qed   
\enddemo   

\vskip13pt

   {\bf Remarks 2.2:}  {\it (1)} \,  Claim  {\it (a)\/}  of Theorem 2.1 
{\sl for  $ M_n $  only\/}  was independently proved in [PW] and in [Ko],
but taking a field as ground ring.  In [Ko], claim  {\it (b)}  {\sl for 
$ {GL}_n $  only\/}  was proved as well.  Similarly, the analogue of
claim  {\it (b)}  {\sl for  $ {SL}_n $  only\/}  was proved in [Ga],
\S 7, but taking as ground ring the field  $ k(q) $   --- for any field 
$ k $  of zero characteristic.  Our proof then provide an alternative,
unifying approach, which yields stronger results over  $ R \, $.
                                                 \par   
    {\it (2)} \,  We would better point out a special aspect of the basic
assumption of Theorem 2.1 about  $ q $  and  $ R \, $.  Namely, if the
subring  $ \langle 1 \rangle $  of  $ R $  generated by  $ 1 $  has
prime characteristic (hence it is a finite field) then the condition
on  $ (q-1) $  is equivalent to  $ q $  being trascendental over  $ R_q $ 
or  $ \, q = 1 \, $.  But if instead the characteristic of  $ \langle
1 \rangle $  is zero or positive non-prime, then  $ (q-1) $  might be
non-invertible in  $ R_q $  even though  $ q $  is algebraic (or even
integral) over  $ \langle 1 \rangle \, $.  
                                            \par   
   The end of the story is that Theorem 2.1 holds true in the ``standard''
case of trascendental values of  $ q \, $,  \, but also in more general
situations.   
                                                 \par   
    {\it (3)} \,  The argument used in the proof of Theorem 2.1 to get
the result for  $ \Oqrsl $  from those for  $ \Oqrgl \, $,  \, via the
isomorphism  $ \, \Oqrsl \otimes_R R\big[x,x^{-1}\big] \cong \Oqrgl \, $, 
\, actually work  {\sl both ways}.  Therefore, one can also prove the
results directly for  $ \Oqrsl $   --- as we sketched above ---   and
from them deduce those for  $ \Oqrgl \, $.  Even more, as we have proved
independently the results for  $ \Oqrgl $   --- i.e.,  $ B^{\,\vee}_{GL} $ 
and  $ B^{\vee,-}_{GL} $  are  $ R $--bases  ---   and  for  $ \Oqrsl $  
--- i.e.,  $ B_{SL} $  and  $ B^-_{SL} $  are  $ R $--bases  ---   we
can use them to prove that the algebra morphism  $ \, \Oqrsl \otimes_R
R\big[x,x^{-1}\big] \longrightarrow \Oqrgl \, $  is in fact bijective.   
                                                 \par   
    {\it (4)} \,  The orders considered in claim  {\it (b)\/}  of Theorem
2.1 refer to a  {\sl triangular decomposition\/}  of  $ \Oqrgl $  and 
$ \Oqrsl $  which is opposite to the standard one.  This opposite
decomposition was introduced   --- and its importance was especially
pointed out ---   in [Ko].   

\vskip17pt   

   We are now ready to state and proof the main result of this paper:   
%
%
 \eject   

\proclaim{Theorem 2.3 (PBW theorem for  $ \Oezg $  as an 
$ \Ozeg $--module,  \text{\rm for}  $ G \! \in \! \big\{
M_n , {GL}_n \big\} $)}  
 \vskip2pt   
   Let any total order be fixed in  $ \, {\{ 1, \dots,
n \}}^{\times 2} \, $.  Then the set of ordered monomials     
 \vskip-5pt   
  $$  \Bb^M_{GL} \,\; := \;\, \bigg\{\; {\textstyle
\prod\limits_{i,j=1}^n} t_{i,j}^{N_{i,j}} \;\bigg\vert\;\,
0 \leq N_{i,j} \leq \ell \! - \! 1 \, , \; \forall \, i,
j \;\bigg\}  $$   
 \vskip-1pt   
\noindent   
 thought of as a subset of  $ \, \Oezm \subset \Oezgl \, $,  \, is a
basis  of  $ \, \Oezm \, $  as a module over  $ \, \Ozem \, $,  \, and
a basis of  $ \, \Oezgl \, $  as a module over  $ \, \Ozegl \, $.   
                                                        \par   
   In particular, both modules are free of rank  $ \ell^{\text{\it dim}(G)} \, $,  \, with  $ \, G \in \big\{ M_n , {GL}_n \big\} \, $.  
\endproclaim   

\demo{Proof}  When specialising,  Theorem 2.1{\it (a)\/}  implies that  $ \Oezm $  is a free  $ \Zeps $--module  with   
 $ \; B_M\Big|_{q=\varepsilon} = \Big\{ \prod_{i,j=1}^n 
t_{ij}^{N_{ij}} \,\Big\vert\; N_{ij} \in \N \;\; \forall \,
i, j \,\Big\} \; $   
as basis   --- where, by abuse of notation, we write again  $ t_{ij} $ 
for  $ \, t_{ij}\big|_{q=\varepsilon} \; $.  Now, whenever the exponent 
$ N_{ij} $  is a multiple of  $ \ell $,  the power  $ t_{ij}^{N_{ij}} $ 
belongs to the isomorphic image  $ \, \gerFr_{\Z} \big( \Ozem \big) \, $ 
of  $ \Ozem $  inside  $ \Oezm $,  hence it is a  {\sl scalar\/}  for
the  $ \Ozem $--module  structure of  $ \Oezm \, $.  Therefore, reducing 
all exponents modulo  $ \ell $  we find that  $ \, \Bb^M_{GL} $  is a spanning set for the  $ \Ozem $--module  $ \Oezm \, $.  In addition, 
$ \Ozm $  clearly admits as  $ \Z $--basis  the set   
 $ \; \overline{B}_M \, = \Big\{ \prod_{i,j=1}^n
\bar{t}_{ij}^{\,N_{ij}} \,\Big\vert\; N_{ij} \in \N \;\;
\forall \, i, j \,\Big\} \; $.  
It follows that  $ \overline{B}_M $  is also a 
$ \Zeps $--basis  of  $ \, \Ozem \, $,  \, so 
 $ \; \gerFr_{\Z} \big(\, \overline{B}_M \big) = \Big\{ \prod_{i,j=1}^n
t_{ij}^{\, \ell \, N_{ij}} \,\Big\vert\; N_{ij} \in \N \;\; \forall
\, i, j \,\Big\} \; $  
 is a  $ \Zeps $--basis  of  $ \, \gerFr_{\Z} \big( \Ozem \big) \, $.  This
last fact easily implies that  $ \, \Bb^M_{GL} \, $  is also  $ \Ozem $--linearly  independent, hence it is a basis of  $ \Oezm $  over  $ \Ozem $  as claimed.   
 \vskip1pt
   As to  $ \Oezgl $,  \, from definitions and the analysis in
\S 1 we get (with  $ \, D_\varepsilon := D_q\big|_\varepsilon \, $)   
 \vskip-9pt   
  $$  \displaylines{ 
   \qquad  \Oezgl  \, = \,  \Oezm\big[D_\varepsilon^{-1}\big]  \, = \, 
\Oezm\big[D_\varepsilon^{-\ell}\,\big]  \, = \,   \hfill  \cr   
   \hfill   \, = \,  \Ozem\big[D^{-1}\big] {\textstyle
\bigotimes\limits_{\Ozem}} \! \Oezm  \, = \,  \Ozegl {\textstyle
\bigotimes\limits_{\Ozem}} \! \Oezm  \qquad  \cr }   $$   
 \vskip-1pt   
\noindent   
 thus the result for  $ \Oezgl $  follows at once from that for 
$ \Oezm \, $.   \qed   
\enddemo   

 \vskip27pt

\centerline { \bf  \S \; 3 \, Frobenius structures }

\vskip13pt

  {\bf 3.1 Frobenius extensions and Nakayama automorphisms.}  Following
[BGStro], we say that a ring  $ R $  is a  {\it free Frobenius
extension\/}  over a subring  $ S $,  if  $ R $  is a free  $ S $--module 
of finite rank, and there is an isomorphism  $ \; F \, \colon R
\longrightarrow \text{Hom}_S(R,S) \, $  of  $ R-S $--bi-modules.  Then 
$ F $  provides a non-degenerate associative  $ S $--bilinear  form  $ \,
\B \, \colon R \times R \longrightarrow S \, $,  \, via  $ \, \B(r,t) =
F(t)(r) \, $.  Conversely, one can characterise Frobenius extensions
using such forms.  When  $ \, S = \Cal{Z} \, $  is contained in the
centre of  $ R \, $,  there is a  $ \Cal{Z} $--algebra  automorphism 
$ \, \nu: R \longrightarrow R \, $,  \, given by  $ \, r \, F(1) = F(1)
\, \nu(r) \, $  (for all  $ \, r \in R \, $),  \, and such  $ \; \B(x,y)
\, = \, \B\big(\nu(y),x\big) \, $.  This is called the  {\it Nakayama
automorphism},  and it is uniquely determined by the pair  $ \, \Cal{Z}
\subseteq R \, $,  \, up to  $ \text{\it Int}\,(R) \, $.   
%
%
%

\vskip13pt

\proclaim{Proposition 3.2}  (cf.~[BGStro], \S 2)
                                                  \par   
   Let  $ R $  be a ring,  $ \Cal{Z} $  an affine central subalgebra
of  $ R \, $.  Assume that  $ R $  is free of finite rank as a 
$ \Cal{Z} $--module,  with a  $ \Cal{Z} $--basis  $ \Cal{B} $  that
satisfies the following condition: there exists a  $ \Cal{Z} $--linear
functional  $ \, \Phi \, \colon R \rightarrow \Cal{Z} \, $  such that
for any non-zero  $ \, a = \sum_{b \in \Cal{B}} z_b b \in R \, $  there
exists  $ \, x \in R \, $  for which  $ \, \Phi(xa) = u z_b \, $  for
some unit  $ \, u \in \Cal{Z} \, $  and some non-zero  $ \, z_b \in
\Cal{Z} \, $.   
                                                  \par   
  Then  $ R $  is a free Frobenius extension of  $ \Cal{Z} $.  Moreover,
for any maximal ideal  $ \frak{m} $  of  $ \Cal{Z} $,  the finite
dimensional quotient  $ R \big/\frak{m} R $  is a finite dimensional
Frobenius algebra.  
\endproclaim
%
%
 \eject   

   This result is used in [BGStro] to show that many families of
algebras   --- in particular, some related to  $ \Oeg \, $,  where 
$ G $  is a (complex, connected, simply-connected) semisimple affine algebraic group ---   are indeed
free Frobenius extensions.  But the authors could not prove the same
for  $ \Oeg \, $,  as they did not know an explicit  $ \Og $--basis 
of  $ \Oeg \, $.  Now, following their strategy  {\sl and\/}  using
Theorem 2.3, I shall now prove that  $ \Oezg $  is free Frobenius
over  $ \Ozeg $  when  $ G $  is  $ M_n $  or  $ {GL}_n \, $.   

\vskip13pt

\proclaim{Theorem 3.3}  Let  $ G $  be  $ M_n $  or  $ {GL}_n \, $.  Then 
$ \Oezg $  is a free Frobenius extension of  $ \, \Ozeg \, $,  \, with
Nakayama automorphism  $ \nu $  given by  $ \; \nu\big(t_{i,j}\big) =
\varepsilon^{2(i+j-n-1)} \, t_{i,j} \; $  ($ \, i, j = 1, \dots, n \, $).     
\endproclaim

\demo{Proof}  We prove that there exists a suitable  $ \Ozeg $--linear 
functional  $ \, \Phi \, \colon \Oezg \longrightarrow \Ozeg \, $  as
required in Proposition 3.2, so that that result applies to  $ \, R
:= \Oezg \, $  and  $ \, \Cal{Z} := \Ozeg \, $.   
                                                  \par   
   Define  $ \Phi $  on the elements of the  $ \Ozeg $--basis 
$ \Bb^M_{GL} $  of  $ \Oezg $  (see Theorem 2.3) by  
  $$  \Phi \bigg(\, {\textstyle \prod\limits_{i,j=1}^n} t_{i,j}^{N_{i,j}}
\bigg)  \, := \;  {\textstyle \prod\limits_{i,j=1}^n} \delta_{N_{i,j},
\ell-1}  \, = \,  \cases   
   1 \, ,  \;\; \text{if}\,\ N_{i,j} = \ell \! - \! 1 
\;\; \forall \; i, j  \\   
   0 \, ,  \;\; \text{if not}   
                  \endcases   \eqno (3.1)  $$    
(for all  $ \; 0 \leq N_{i,j} \leq \ell \! - \! 1 \, $),  and extend to 
all of  $ \Oezg $  by  $ \Ozeg $--linearity.  In other words,  $ \Phi $ 
is the unique  $ \Ozeg $--valued  linear functional on  $ \Oezg $  whose
value is 1 on the basis element  $ \, \t^{\,\underline{\ell-1}} :=
\prod\limits_{i,j=1}^n t_{i,j}^{\,\ell-1} \, $  and is zero on
all other elements of the  $ \Ozeg $--basis  $ \Bb^M_{GL} \, $.   
 \vskip1pt   
   We claim that  $ \Phi $  satisfies the assumptions of Proposition 3.2, so
the latter applies and proves our statement.  Indeed, let us consider any
non-zero  $ \, a = \sum_{\t \in \Bb^M_{GL}} z_{\t} \, \t \in \Oezg \, $,  \,
and let  $ \, \t_{\,0} = \prod\limits_{i,j=1}^n t_{i,j}^{N_{i,j}} \, $  in 
$ \Bb^M_{GL} $  be such that  $ \, z_{\t_{\,0}} \not= 0 \, $  and  $ w(\t_{\,0}) $ 
is maximal (w.r.t.~$ \leq_{lex} $).  Then define  $ \; \t_{\,0}^\vee
:= \, \prod\limits_{i,j=1}^n t_{i,j}^{N'_{i,j}} \; $  $ \Big(\! \in
\B^M_{GL} \Big) $  with  $ \, N'_{i,j} := \ell - 1 - N_{i,j} \, $ 
for all  $ i, j = 1, \dots, n \, $.  Quoting from the proof of 
Theorem 2.1{\it (a)},  we know that  $ \; \t_{\,0}^\vee \, \t_{\,0} \,
= \, \varepsilon^s \, \t^{\,\underline{\ell-1}} \, + \, \text{\it l.t.'s}
\, $,  \; where  $ \, s \in \Z \, $  and the expression  {\it l.t.'s\/} 
now stands for an  $ \Ozeg $--linear  combination of monomials  $ \,
\underline{\check{t}} \in \Bb^M_{GL} \, $  such that  $ \, w \big(
\underline{\check{t}} \,\big) \lneqq_{lex} w \big( \t^{\,
\underline{\ell-1}} \big) \, $;  \, in particular,  $ \,
\Phi \big( \underline{\check{t}} \,\big) = 0 \, $  for all
these  $ \underline{\check{t}} \, $,  \, hence eventually  $ \,
\Phi \big( \t_{\,0}^\vee \, \t_{\,0} \big) = \varepsilon^s \, \Phi
\big( \t^{\,\underline{\ell-1}} \big) = \varepsilon^s \, $.  Similarly,
if  $ \, \t' \in \B^M_{GL} \, $  is such that  $ \, w\big(\t'\big)
<_{lex} w(\t) \, $,  \, then  $ \, \t_{\,0}^\vee \t' \, $  is an 
$ \Ozeg $--linear  combination of PBW monomials whose weight
is at most  $ w\big(\t_{\,0}^\vee \, \t'\big) $,  hence  $ \, \Phi
\big( \t_{\,0}^\vee \, \t' \big) = 0 \, $.  As we chose  $ \t_{\,0} $ 
so that  $ w(\t_{\,0}) $  is maximal, we eventually find   
  $$  \Phi \big( \t_{\,0}^\vee \, a \big)  \; = \;  {\textstyle
\sum_{\t \in \Bb^M_{GL}}} z_{\t} \, \Phi(\t)  \; = \;  z_{\t_{\,0}}
\, \Phi(\t_{\,0})  \; = \;  \varepsilon^s z_{\t_{\,0}}  $$   
where  $ \varepsilon^s $  is a unit in  $ \Ozeg \, $.  So  $ \Phi $ 
satisfies the assumptions of Proposition 3.2, as claimed.
 \vskip3pt   
   As to the Nakayama automorphism  $ \, \nu \, \colon \Oezg
\longrightarrow \Oezg \, $,  \, it is characterized (see \S 3.1)
by the property that  $ \; \B(x,y) \, = \, \B\big(\nu(y),x\big)
\; $  for all  $ \, x $, $ y \in R \, $.  Here  $ \B $  is a 
$ \Cal{Z} $--bilinear  form as in \S 3.1, which now is related
to  $ \Phi $  by the formula  $ \; \B(x,y) = \Phi (x y) \; $ 
for all  $ x, y \in R \, $.   
                                               \par   
   As  $ \Phi $  is an automorphism, and  $ \Oezg $  is generated  
--- over  $ \Ozeg $  ---   by the  $ t_{i,j} $'s,  the claim about 
$ \nu $  is proved if we show that   
  $$  \Phi \Big( {\textstyle \prod_{r,s=1}^n} t_{r,s}^{e_{r,s}}
\cdot t_{i,j} \Big)  \; = \;  \Phi \Big( \varepsilon^{2(i+j-n-1)}
\, t_{i,j} \cdot {\textstyle \prod_{r,s=1}^n} t_{r,s}^{e_{r,s}} \Big)  
\eqno (3.2)  $$   
   \indent   Now, our usual argument shows that the expansions of the
product of a generator  $ t_{i,j} $  and a PBW monomial  $ \prod_{r,s=1}^n
t_{r,s}^{e_{r,s}} $  (in either order of the factors) as an 
$ \Ozeg $--linear  combination of elements of the  $ \Ozeg $--basis 
$ \Bb^M_{GL} $  are of the form   
 \vskip-4pt  
  $$  \eqalign{ 
   {\textstyle \prod\limits_{r,s=1}^n} t_{r,s}^{e_{r,s}}
\cdot t_{i,j}  &  \; = \;  \varepsilon^{i+j-2n} \, {\textstyle
\prod\limits_{r,s=1}^n}
t_{r,s}^{e_{r,s} + \delta_{r,i} \delta_{j,s}} \, + \, \text{\it l.t.'s} 
\cr  
   t_{i,j} \cdot {\textstyle \prod\limits_{r,s=1}^n} t_{r,s}^{e_{r,s}} 
&  \; = \;  \varepsilon^{2-i-j} \, {\textstyle \prod\limits_{r,s=1}^n}
t_{r,s}^{e_{r,s} + \delta_{r,i} \delta_{j,s}} \, + \, \text{\it l.t.'s} 
\cr }  $$   
 \vskip-4pt  
\noindent   
 This along with (3.1) gives   
  $$  \displaylines{ 
   \Phi \bigg(\, {\textstyle \prod\limits_{r,s=1}^n} t_{r,s}^{e_{r,s}}
\cdot t_{i,j} \bigg)  = \,  \varepsilon^{i+j-2n} \, \Phi \bigg(\,
{\textstyle \prod\limits_{r,s=1}^n} t_{r,s}^{e_{r,s} + \delta_{r,i}
\delta_{j,s}} \bigg)  = \,  \varepsilon^{i+j-2n}  \quad  \text{if\ } 
\; e_{r,s} = \ell - 1 - \delta_{r,i} \, \delta_{j,s}  \cr 
   \Phi \bigg(\, {\textstyle \prod\limits_{r,s=1}^n} t_{r,s}^{e_{r,s}}
\cdot t_{i,j} \bigg)  = \,  \varepsilon^{i+j-2n} \, \Phi \bigg(\,
{\textstyle \prod\limits_{r,s=1}^n} t_{r,s}^{e_{r,s} + \delta_{r,i}
\, \delta_{j,s}} \bigg)  = \,  0  \hskip37pt  \text{if not} 
\hskip81pt  \cr }  $$   
and similarly    
  $$  \displaylines{ 
   \Phi \bigg(\, t_{i,j} \cdot {\textstyle \prod\limits_{r,s=1}^n}
t_{r,s}^{e_{r,s}} \bigg)  = \,  \varepsilon^{2-i-j} \, \Phi \bigg(\,
{\textstyle \prod\limits_{r,s=1}^n} t_{r,s}^{e_{r,s} + \delta_{r,i}
\delta_{j,s}} \bigg)  = \,  \varepsilon^{2-i-j}  \quad  \text{if\ } 
\; e_{r,s} = \ell - 1 - \delta_{r,i} \, \delta_{j,s}  \cr 
   \Phi \bigg(\, t_{i,j} \cdot {\textstyle \prod\limits_{r,s=1}^n}
t_{r,s}^{e_{r,s}} \bigg)  = \,  \varepsilon^{2-i-j} \, \Phi \bigg(\,
{\textstyle \prod\limits_{r,s=1}^n} t_{r,s}^{e_{r,s} + \delta_{r,i}
\, \delta_{j,s}} \bigg)  = \,  0  \hskip32pt  \text{if not}  
\hskip80pt  \cr }  $$   
Direct comparison now shows that (3.2) holds, q.e.d.   \qed   
\enddemo  

 \vskip21pt

\centerline{ \sc acknowledgements }

\vskip1pt

\centerline{ \smallrm The author thanks I.~Gordon, Z.~Skoda,
  C.~Stroppel and the referee for several useful comments. }

 \vskip1,1truecm

\vskip9pt

\enddocument
\bye
\end